\documentclass{amsart}
\usepackage[all]{xy}
\usepackage{color}
\usepackage{amsthm}
\usepackage{amssymb}
\usepackage[colorlinks=true]{hyperref}
\numberwithin{equation}{section}

\newtheorem{theorem}[equation]{Theorem}
\newtheorem*{theorem*}{Theorem}
\newtheorem{lemma}[equation]{Lemma}

\newtheorem*{conjecture*}{Mamma Conjecture}
\newtheorem*{conjecture1*}{Mamma Conjecture (revisited)}
\newtheorem{proposition}[equation]{Proposition}

\newtheorem*{corollary*}{Corollary}

\theoremstyle{remark}
\newtheorem{definition}[equation]{Definition}

\newtheorem{example}[equation]{Example}

\newtheorem{notation}[equation]{Notation}

\theoremstyle{remark}

\newcommand{\cA}{{\mathcal A}}

\newcommand{\cC}{{\mathcal C}}

\newcommand{\cN}{{\mathcal N}}
\newcommand{\cO}{{\mathcal O}}

\newcommand{\cR}{{\mathcal R}}

\newcommand{\sVect}{\mathrm{sVect}}
\newcommand{\Vect}{\mathrm{Vect}}

\newcommand{\bbG}{\mathbb{G}}

\newcommand{\bbZ}{\mathbb{Z}}

\DeclareMathOperator{\NChow}{NChow} % category of noncommutative Chow motives
\DeclareMathOperator{\NNum}{NNum} % category of noncommutative numerical motives
\DeclareMathOperator{\NHom}{NHom} % category of noncommutative homological motives
\DeclareMathOperator{\Chow}{Chow} % category of Chow motives
\DeclareMathOperator{\Num}{Num} % category of numerical motives
\DeclareMathOperator{\even}{even} % category of numerical motives
\DeclareMathOperator{\odd}{odd} % category of numerical motives

\DeclareMathOperator{\Gal}{Gal} % Galois group

\newcommand{\perf}{\mathrm{perf}}

\newcommand{\dg}{\mathrm{dg}}

\newcommand{\Hom}{\mathrm{Hom}}

\newcommand{\too}{\longrightarrow}

\newcommand{\ie}{\textsl{i.e.}\ }

\title[Unconditional noncommutative motivic Galois groups]{Unconditional\\ noncommutative motivic Galois groups}
\author{Matilde Marcolli and Gon{\c c}alo~Tabuada}

\address{Matilde Marcolli, Mathematics Department, Mail Code 253-37, Caltech, 1200 E.~California Blvd. Pasadena, CA 91125, USA}
\email{matilde@caltech.edu} 
\urladdr{http://www.its.caltech.edu/~matilde}

\address{Gon{\c c}alo Tabuada, Department of Mathematics, MIT, Cambridge, MA 02139, USA}
\email{tabuada@math.mit.edu}
\urladdr{http://math.mit.edu/~tabuada}
\thanks{M.~Marcolli was partially supported by NSF grants DMS-0901221, DMS-1007207, DMS-1201512, and PHY-1205440. G.~Tabuada was partially supported by the NEC Award-2742738.} 
\thanks{The authors are very grateful to the organizers of the conference ``Hodge Theory and Classical Algebraic Geometry'' 
for the opportunity to present this work.}
\date{\today}


\begin{document}
\begin{abstract}
In this short note we introduce the unconditional noncommutative motivic Galois groups and relate them with those of Andr{\'e}-Kahn. 
\end{abstract}

\maketitle
\vskip-\baselineskip
%\vskip-\baselineskip
%\vskip-\baselineskip
%\tableofcontents
\smallskip

{\em Dedicated to Herb Clemens, on the occasion of his non-retirement.}
%--------------------------------------------------------
\section*{Motivating questions}\label{Intro}
%--------------------------------------------------------
Motivic Galois groups were introduced by Grothendieck in the sixties as part of his broad theory of (pure) motives. These group schemes are {\em conditional} in the sense that their construction makes use of the standard conjectures. Thanks to Kontsevich \cite{IAS,Miami,finMot}, Grothendieck's theory of motives admits a noncommutative counterpart, with schemes replaced by dg categories. The standard conjectures admit noncommutative analogues and there exist also conditional noncommutative motivic Galois groups; consult \cite{Konts,Semi,Galois} for details. 

Recently, via an evolved ``$\otimes$-categorification'' of the Wedderburn-Malcev's theorem, Andr{\'e}-Kahn \cite{AK-cras} introduced  {\em unconditional} motivic Galois groups. These group schemes $\Gal_H$ (attached to a Weil cohomology theory $H$) are well-defined up to an interior automorphism and do not require the assumption of any of Grothendieck's standard conjectures. This lead us naturally to the following motivating questions:

\smallbreak
{\bf Q1:} \textit{Do the Galois groups $\Gal_H$ admit noncommutative analogues $\Gal_H^{NC}$~?}

{\bf Q2:} \textit{What is the relation between $\Gal_H$ and $\Gal_H^{NC}$~?}

{\bf Q3:} \textit{What is the relation between $\Gal_H^{NC}$ and the conditional Galois groups ?}

\smallbreak

In this short note we provide precise answers to these three questions; see Definition \ref{def:unconditional}, Theorem~\ref{thm:main1}, and Proposition~\ref{prop:very-last}, respectively. 

%--------------------------------------------------------
\section{Preliminaries}\label{sec:preliminaries}
%--------------------------------------------------------
Let $k$ be a base field and $F$ a field of coefficients. The classical idempotent completion construction will be denoted by $(-)^\natural$.
%--------------------------------------------------------
\subsection*{Motives}
%--------------------------------------------------------
We assume the reader is familiar with the categories of Chow motives $\Chow(k)_F$, homological motives $\Hom(k)_F$, and numerical motives $\Num(k)_F$; consult \cite[\S4]{Andre}. The Tate motive will be denoted by $F(1)$. At \S\ref{sec:conditional} we will assume some familiarity with Grothendieck's standard conjectures of type $C$ (K\"unneth) and $D$ (homological=numerical) as well as with the sign conjecture $C^+$; see \cite[\S5]{Andre}. 
%--------------------------------------------------------
\subsection*{Dg categories}
%--------------------------------------------------------
A differential graded (=dg) category $\cA$ is a category enriched over cochain complexes of $k$-vector spaces; consult \cite{ICM} for details. Every (dg) $k$-algebra $A$ gives naturally rise to a dg category $\underline{A}$ with a single object. Another source of examples is provided by schemes since the category of perfect complexes $\perf(X)$ of every quasi-compact quasi-separated $k$-scheme $X$ admits a canonical dg enhancement $\perf_\dg(X)$; see \cite[\S4.6]{ICM}. When $X$ is quasi-projective this dg enhancement is moreover unique; see \cite[Thm.~2.12]{LO}. 
%--------------------------------------------------------
\subsection*{Noncommutative motives}
%--------------------------------------------------------
We assume the reader is familiar with the categories of noncommutative Chow motives $\NChow(k)_F$, noncommutative homological motives $\NHom(k)_F$, and noncommutative numerical motives $\NNum(k)_F$; consult the survey articles \cite[\S 2-3]{survey} \cite[\S4]{survey1} and the references therein. At \S\ref{sec:conditional} we will assume some familiarity with conjectures $C_{NC}$ (the noncommutative analogue of $C^+$) and $D_{NC}$ (the noncommutative analogue of type $D$); see \cite[\S4]{survey}.
%-----------------------------------------------------------------------
\section{Construction of the unconditional NC motivic Galois groups}\label{sec:Galois}
%-----------------------------------------------------------------------
Assume that $k$ is of characteristic zero and that $k \subseteq F$ or $F \subseteq k$. As proved in \cite[Thm.~7.2]{Galois}, periodic cyclic homology $HP$ gives rise to an $F$-linear $\otimes$-functor
\begin{equation}\label{eq:HP}
HP_\ast:\NChow(k)_F \too \sVect(K)
\end{equation}
with values in the category of finite dimensional super $K$-vector spaces (with $K=F$ when $k \subseteq F$ and $K=k$ when $F \subseteq k$). By definition of the category of noncommutative homological motives, \eqref{eq:HP} descends to a faithful $F$-linear $\otimes$-functor
\begin{equation}\label{eq:HP1}
HP_\ast: \NHom(k)_F \too \sVect(K)\,.
\end{equation}
\begin{notation}
Let $\NHom(k)_F^\pm$ be the full subcategory of those noncommutative homological motives $N$ whose associated K{\"u}nneth projectors
$$ \pi^\pm_N: HP_\ast(N) \twoheadrightarrow HP_\ast^\pm(N) \hookrightarrow HP_\ast(N)$$
can be written as $\pi^\pm_N=HP_\ast(\underline{\pi}^\pm_N)$ with $\underline{\pi}^\pm_N$ endomorphisms in $\NHom(k)_F$. 
\end{notation}
\begin{example}
Let $X$ be a smooth projective $k$-scheme. When $F \subseteq k$, the proof of \cite[Thm.~1.3]{Galois} shows us that $\perf_\dg(X)$ belongs to $\NHom(k)^\pm_F$ whenever $X$ satisfies the sign conjecture. As proved by Kleiman~\cite{Kleiman}, this holds when $X$ is an abelian variety. Using the stability of $\NHom(k)^\pm_F$ under direct factors and tensor products (see \cite[Prop.~8.2]{Galois}) we then obtain a large class of examples.
\end{example}
Since by hypothesis $k$ (and hence $F$) is of characteristic zero, \cite[Prop.~2]{AK-cras} applied to the above functor \eqref{eq:HP1} gives rise to a new rigid symmetric monoidal category $\NHom^\dagger(k)_F^\pm$ (obtained from $\NHom(k)_F^\pm$ by modifying its symmetric isomorphism constraints) and to a composed faithful $F$-linear $\otimes$-functor
\begin{equation}\label{eq:HP2}
\NHom^\dagger(k)_F^\pm \subset \NHom(k)_F \stackrel{HP_\ast}{\too} \sVect(K) \stackrel{\mathrm{forget}}{\too} \Vect(K)\,.
\end{equation}
More importantly, $\NHom^\dagger(k)_F^\pm$ is semi-primary, its $\otimes$-ideals $\cR^\pm$ and $\cN^\pm$ agree, the quotient category $\NNum^\dagger(k)^\pm_F$ defined by $\NHom^\dagger(k)^\pm_F/\cN^\pm $ is abelian semi-simple, and the canonical projection $\otimes$-functor
\begin{equation}\label{eq:proj}
\NHom^\dagger(k)^\pm_F \too \NNum^\dagger(k)^\pm_F
\end{equation}
is conservative; consult \cite[\S1]{AK-cras}. Thanks to \cite[Thm.~8 a)]{AK-cras}, the projection \eqref{eq:proj} admits a $\otimes$-section and any two such $\otimes$-sections are conjugated by a $\otimes$-isomorphism. The choice of a $\otimes$-section $s^{NC}$ gives then rise to a faithful $F$-linear $\otimes$-functor
\begin{equation}\label{eq:fiber-last}
f_{HP}: \NNum^\dagger(k)_F^\pm \stackrel{s^{NC}}{\too} \NHom^\dagger(k)_F^\pm \stackrel{\eqref{eq:HP2}}{\too} \Vect(K)\,.
\end{equation}
The category $\NNum^\dagger(k)^\pm_F$, endowed with the fiber functor $f_{HP}$, becomes a Tannakian category; see \cite[Appendix~A]{Galois}. 
\begin{definition}\label{def:unconditional}
The {\em unconditional noncommutative motivic Galois group} $\Gal_{HP}^{NC}$ is the group scheme $\underline{\mathrm{Aut}}^\otimes(f_{HP})$ of $\otimes$-automorphisms of the above fiber functor \eqref{eq:fiber-last}.
\end{definition}

A different choice of the $\otimes$-section $s^{NC}$ gives rise to an isomorphic group scheme (via an interior automorphism). Moreover, since $\NHom^\dagger(k)^\pm_F$ is abelian semi-simple, the group scheme $\Gal^{NC}_{HP}$ is {\em pro-reductive}, \ie its unipotent radical is trivial.
%-----------------------------------------------------------------------
\section{Relation with Andr{\'e}-Kahn's motivic Galois groups}\label{sec:AK}
%-----------------------------------------------------------------------
Assume that $k$ is of characteristic zero and that $F \subseteq k$. As explained by Andr{\'e}-Kahn \cite{AK-cras}, de Rham cohomology theory $H^\ast_{dR}$ gives rise to a well-defined fiber functor
$$f_{dR}: \Num^\dagger(k)^\pm_F \too \mathrm{Vect}(k)$$
and consequently to an unconditional motivic Galois group $\Gal_{dR}:=\underline{\mathrm{Aut}}^\otimes(f_{dR})$. Let us denote by $\langle F(1)\rangle$ the Tannakian subcategory of $\Num^\dagger(k)^\pm_F$ generated by the Tate motive $F(1)$ and write $\Gal_{dR}(F(1))$ for the group scheme of $\otimes$-automorphisms of the composed fiber functor
\begin{equation}\label{eq:fiber-Tate}
\langle F(1) \rangle \stackrel{t}{\hookrightarrow} \Num^\dagger(k)^\pm_F \stackrel{f_{dR}}{\too} \Vect(k)\,.
\end{equation}
As explained in \cite[\S2.3.3]{Andre}, the inclusion of categories gives rise to an homomorphism $t: \Gal_{dR} \twoheadrightarrow \Gal_{dR}(F(1))$. The relation between $\Gal_{dR}$ and $\Gal_{HP}^{NC}$ is the following:
\begin{theorem}\label{thm:main1}
There exists a comparison group scheme homomorphism
\begin{equation}\label{eq:comparison}
\Gal_{HP}^{NC} \to \mathrm{Kernel}(t: \Gal_{dR} \twoheadrightarrow \Gal_{dR}(F(1)))\,.
\end{equation}
In the case where $k=F$, $\Gal_{dR}(F(1))$ identifies with the multiplicative group scheme $\mathbb{G}_m$ and the above comparison homomorphism \eqref{eq:comparison} is faithfully flat
\begin{equation}\label{eq:comparison1}
\Gal_{HP}^{NC} \twoheadrightarrow \mathrm{Kernel}(t: \Gal_{dR} \twoheadrightarrow \mathbb{G}_m)\,.
\end{equation}
\end{theorem}
Intuitively speaking, \eqref{eq:comparison1} shows us that the $\otimes$-symmetries of the commutative world which can be lifted to the noncommutative world are precisely those that become trivial when restricted to the Tate motive. It is unclear at the moment if the kernel of the comparison homomorphisms \eqref{eq:comparison}-\eqref{eq:comparison1} is non-trivial.
\begin{proof}
We start by constructing the comparison homomorphism \eqref{eq:comparison}. Recall from the proof of \cite[Thm.~1.7]{Galois} that we have the following commutative diagram
\begin{equation}\label{eq:diag-global}
\xymatrix@C=1.5em@R=2em{
\Chow(k)_F \ar[r] \ar[d] & \Chow(k)_F\!/_{\!\!-\otimes F(1)} \ar[d] \ar[r] & \NChow(k)_F \ar[d] \ar[r]^-{HP_\ast} & \mathrm{sVect}(k) \ar@{=}[d]\\
\Chow(k)_F\!/\mathrm{Ker} \ar[d] \ar[r] & (\Chow(k)_F\!/_{\!\!-\otimes F(1)})\!/\mathrm{Ker} \ar[r] \ar[d] & \NHom(k)_F \ar[r]_{HP_\ast} \ar[d] & \mathrm{sVect}(k)\\
\Num(k)_F \ar[r] & \Num(k)_F\!/_{\!\!-\otimes F(1)} \ar[r] & \NNum(k)_F& \,,
}
\end{equation}
where $\mathrm{Ker}$ stands for the kernel of the respective horizontal composition towards $\sVect(k)$. Since by hypothesis $k$ is of characteristic zero, the proof of \cite[Thm.~1.3]{Galois} shows us that the upper horizontal composition in \eqref{eq:diag-global} identifies with the functor
\begin{eqnarray}\label{eq:deRham}
&sH^\ast_{dR} : \Chow(k)_F \too \sVect(k) & X \mapsto (\underset{n \even}\oplus H^n_{dR}(X), \underset{n \odd}\oplus H^n_{dR}(X))\,.
\end{eqnarray}
Hence, its kernel $\mathrm{Ker}$ agrees with the one of de Rham cohomology theory
\begin{eqnarray*}
H^\ast_{dR}: \Chow(k)_F \too \mathrm{GrVect}(k)_{\geq 0} & &X \mapsto \{H^n_{dR}(X)\}_{n \geq 0}\,.
\end{eqnarray*}
 As a consequence, the idempotent completion of $\Chow(k)_F\!/\mathrm{Ker}$ agrees with the category $\Hom(k)_F$; see \cite[\S4]{Andre}. The two lower commutative squares in \eqref{eq:diag-global}, combined with the fact that $\Num(k)_F$ and $\NHom(k)_F$ are idempotent complete, give then rise to the following commutative diagram
\begin{equation}\label{eq:com-radical}
\xymatrix{
\Hom(k)_F \ar[d] \ar[r]^-\Phi & \NHom(k)_F \ar[d] \\
\Num(k)_F \ar[r] & \NNum(k)_F \,. 
} 
 \end{equation}
By construction of $\Num(k)_F$ and $\NNum(k)_F$, the kernels of the vertical functors in \eqref{eq:com-radical} are precisely the largest $\otimes$-ideals of $\Hom(k)_F$ and $\NHom(k)_F$. Hence, the commutativity of diagram \eqref{eq:com-radical} allows us to conclude that the functor $\Phi$ is {\em radical}, \ie that it preserves these largest $\otimes$-ideals. 

Now, note that by construction, the following composition
$$ \Hom(k)_F \stackrel{\Phi}{\too} \NHom(k)_F \stackrel{HP_\ast}{\too} \sVect(k)$$
agrees with the factorization of the above functor \eqref{eq:deRham} through the category $\Hom(k)_F$. In particular, it is faithful. By applying \cite[Prop.~2]{AK-cras} to $HP_\ast$ and $HP_\ast \circ \Phi$, we obtain then an induced $F$-linear $\otimes$-functor $\Phi: \Hom^\dagger(k)^\pm_F \to \NHom^\dagger(k)^\pm_F$. This functor is also radical and therefore gives rise to the following solid commutative square (the dotted arrows denote the $\otimes$-sections provided by \cite[Thm.~8 a)]{AK-cras})
\begin{equation*}%\label{eq:commutative1}
\xymatrix{
\Hom^\dagger(k)_F^\pm \ar[d] \ar[rr]^-\Phi && \NHom^\dagger(k)_F^\pm \ar[d] \\
\Num^\dagger(k)_F^\pm \ar[rr]_-{\overline{\Phi}} \ar@{.>}@/^1.5pc/[u]^s && \ar@{.>}@/_1.5pc/[u]_{s^{NC}} \NNum^\dagger(k)_F^\pm\,.
}
\end{equation*}
Note that $\Gal_{dR}$ can be described as the $\otimes$-automorphisms of the fiber functor
$$f'_{dR}: \Num^\dagger(k)^\pm_F \stackrel{s}{\too} \Hom^\dagger(k)^\pm_F \stackrel{\Phi}{\too} \NHom^\dagger(k)^\pm_F \stackrel{\eqref{eq:HP2}}{\too} \Vect(k)\,.$$
Since $\Phi$ is radical, \cite[Props.~12.2.1 and 13.7.1]{AK} imply that the two $\otimes$-functors
\begin{eqnarray*}
 \Phi \circ s:  \Num^\dagger(k)_F^\pm \too \NHom^\dagger(k)_F^\pm
 && s^{NC} \circ \overline{\Phi}: \Num^\dagger(k)_F^\pm \too \NHom^\dagger(k)_F^\pm
  \end{eqnarray*}
are naturally isomorphic (via a $\otimes$-isomorphism). As a consequence, the above fiber functor $f_{dR}$ becomes naturally $\otimes$-isomorphic to the following composition
$$ \Num^\dagger(k)^\pm_F \stackrel{\overline{\Phi}}{\too} \NNum^\dagger(k)^\pm_F \stackrel{f_{HP}}{\too} \Vect(k)\,.$$
Hence, by definition of the unconditional motivic Galois groups, the functor $\overline{\Phi}$ gives rise to a well-defined comparison group scheme homomorphism $\Gal^{NC}_{HP} \to \Gal_{dR}$. It remains then to show that the composition
\begin{equation}\label{eq:composition}
\Gal_{H}^{NC} \too \Gal_{dR} \stackrel{t}{\twoheadrightarrow} \Gal_{dR}(F(1))
\end{equation}
is trivial. Since the categories $\Num^\dagger(k)^\pm_F$ and $\NNum^\dagger(k)^\pm_F$ are abelian semi-simple and $\overline{\Phi}$ is $F$-linear and additive, we conclude that $\overline{\Phi}$ is moreover {\em exact}, \ie that it preserves kernels and cokernels. Thanks to the commutative diagram \eqref{eq:diag-global} we observe that the image of $F(1)$ under $\overline{\Phi}$ is precisely the $\otimes$-unit $\underline{k}$ of $\Num^\dagger(k)^\pm_F$. Hence, since $\Phi$ is also symmetric monoidal, we obtain the commutative diagram
\begin{equation}\label{eq:com2}
\xymatrix{
*+<2.5ex>{\langle F(1) \rangle} \ar[r]^-{\overline{\Phi}} \ar@{_{(}->}[d]^-t & *+<2.5ex>{\langle \underline{k} \rangle} \ar@{^{(}->}[d] \\
\Num^\dagger(k)^\pm_F \ar[r]_-{\overline{\Phi}}& \NNum^\dagger(k)^\pm_F\,,
}
\end{equation}
where $\langle \underline{k}\rangle$ denotes the Tannakian subcategory of $\NNum^\dagger(k)^\pm_F$ generated by $\underline{k}$. The group scheme of $\otimes$-automorphisms of the fiber functor
$$ \langle \underline{k} \rangle \hookrightarrow \NNum^\dagger(k)^\pm_F \stackrel{f_{HP}}{\too} \Vect(k)$$
is clearly trivial. Using the commutativity of diagram \eqref{eq:com2}, we conclude finally that the above composition \eqref{eq:composition} is trivial. 

Let us now assume that $k=F$. Note first that $(\Num^\dagger(k)^\pm_F, w, F(1))$ is a Tate sub-triple of the one described in \cite[Example~A.5(i)]{Galois} ($w$ stands for the weight $\bbZ$-grading). Since by hypothesis $k=F$ this Tate triple is moreover neutral, with fiber functor given by $f'_{dR}$. As explained in the proof of \cite[Prop.~11.1]{Galois}, the weight $\bbZ$-grading of the Tate triple furnish us the following factorization
$$
\xymatrix{
\langle F(1) \rangle \ar[r]^-{\overline{\eqref{eq:fiber-Tate}}} \ar[dr]_-{\eqref{eq:fiber-Tate}} & \mathrm{GrVect}(k)_\bbZ \ar[d]^-{\mathrm{forget}} \\
& \Vect(k)\,,
}
$$
where $\mathrm{GrVect}(k)_\bbZ$ stands for the category of finite dimensional $\bbZ$-graded $k$-vector spaces. An argument similar to the one in {\em loc.~cit.} allows us then to conclude that $\Gal_{dR}(F(1)) \simeq \bbG_m$. Let us now prove that the comparison homomorphism \eqref{eq:comparison1} is faithfully flat. By applying the general \cite[Prop.~11.1]{Galois} to the neutral Tate-triple $(\Num^\dagger(k)^\pm_F, w, F(1))$ we obtain the following group scheme isomorphism
\begin{equation}\label{eq:isom1}
\Gal((\Num^\dagger(k)^\pm_F\!/_{\!\!-\otimes F(1)})^\natural) \stackrel{\sim}{\too} \mathrm{Kernel}(t: \Gal_{dR} \twoheadrightarrow \bbG_m)\,,
\end{equation} 
where the left-hand-side is the group scheme of $\otimes$-automorphisms of the functor
$$ (\Num^\dagger(k)^\pm_F\!/_{\!\!-\otimes F(1)})^\natural \stackrel{\Psi}{\too} \NNum^\dagger(k)^\pm_F \stackrel{f_{HP}}{\too} \Vect(k)\,.$$
We claim that $\Psi$ is fully-faithful. Consider the composition at the bottom of diagram \eqref{eq:diag-global}, and recall that the functor $\Num(k)_F /_{\!\!-\otimes F(1)} \to \NNum(k)_F$ is fully-faithful. By first restricting ourselves to $\Num(k)^\pm_F$ and then by modifying the symmetry isomorphism constraints we obtain the following composition
\begin{equation*}
\Num^\dagger(k)^\pm_F \too (\Num(k)^\pm_F\!/_{\!\!-\otimes F(1)})^\dagger \stackrel{\psi}{\too} \NNum^\dagger(k)^\pm_F\,.
\end{equation*}
The general \cite[Lemma~B.9]{Galois} (applied to $\cC=\Num(k)^\pm_k$ and $\cO=F(1)$) furnish us a canonical $\otimes$-equivalence between $\Num^\dagger(k)^\pm_F\!/_{\!\!-\otimes F(1)}$ and $(\Num(k)_F^\pm\!/_{\!\!-\otimes F(1)})^\dagger$. Hence, since $\psi$ is fully-faithful and $\NNum^\dagger(k)^\pm_F$ is idempotent complete we conclude that $\Psi$ is also fully-faithful. As a consequence, we obtain a faithfully flat homomorphism
$$ \Gal_{HP}^{NC} \twoheadrightarrow \Gal((\Num^\dagger(k)^\pm_F\!/_{\!\!-\otimes F(1)})^\natural ) \,.$$
By combining it with \eqref{eq:isom1} we conclude finally that the comparison homomorphism \eqref{eq:comparison} is faithfully flat. This achieves the proof.
\end{proof}
%-----------------------------------------------------------------------
\section{Relation with the conditional motivic Galois groups.}\label{sec:conditional}
%-----------------------------------------------------------------------
Assuming the standard conjectures of type $C$ and $D$, we have well-defined conditional motivic Galois groups $\Gal(\Num^\dagger(k)_F)$; see \cite[\S6]{Andre}. Similarly, assuming conjectures $C_{NC}, D_{NC}$, we have well-defined conditional noncommutative motivic Galois groups $\Gal(\NNum^\dagger(k)_F)$; see \cite[\S6]{survey}. As explained in \cite{AK-cras}, $\Gal(\Num^\dagger(k)_F)$ identifies with $\Gal_{dR}$. In the noncommutative world, the following holds:
\begin{lemma}
The conditional Galois group $\Gal(\NNum^\dagger(k)_F)$ identifies with $\Gal_{HP}^{NC}$.
\end{lemma}
\begin{proof}
Conjecture $C_{NC}$ implies that $\NNum^\dagger(k)_F = \NNum^\dagger(k)^\pm_F$. On the other hand, conjecture $D_{NC}$ allows us to choose for the $\otimes$-section $s^{NC}$ the identity functor. This achieves the proof.
\end{proof}
In what concerns the comparison homomorphism \eqref{eq:comparison1}, the following holds:
\begin{proposition}\label{prop:very-last}
The conditional comparison homomorphism
\begin{equation}\label{eq:ind-conj}
\Gal(\NNum^\dagger(k)_F)  \twoheadrightarrow \mathrm{Kernel}(t: \Gal(\Num^\dagger(k)_F)  \twoheadrightarrow \mathbb{G}_m)\,,
\end{equation}
constructed originally in \cite[Thm.~1.7]{Galois}, identifies with \eqref{eq:comparison1}.
\end{proposition}
\begin{proof}
It follows from the fact that the comparison homomorphisms \eqref{eq:ind-conj} and \eqref{eq:comparison1} are induced by the same functor from numerical motives to noncommutative numerical motives, namely by the bottom horizontal composition in diagram \eqref{eq:diag-global}. 
\end{proof}
\smallbreak
\noindent\textbf{Acknowledgments:} The authors are grateful to Michael Artin, Eric Friedlander, Steven Kleiman, and Yuri Manin for useful discussions and motivating questions.

\end{document}